\documentclass[a4paper,10pt,leqno]{article}
\usepackage{amsfonts}
\usepackage{amsmath, amsfonts, amssymb, amsthm}
\numberwithin{equation}{section}

\begin{document}

\title{Lipschitz metric for the modified two-component Camassa-Holm system}

\author
{Chunxia $\mbox{Guan}^{1}$ \footnote{E-mail: guanchunxia123@163.com},
\quad Kai  $\mbox{Yan}^{2}$ \footnote{E-mail: kaiyan@hust.edu.cn, corresponding author}
\quad and \quad
Xuemei $\mbox{Wei}^{1}$ \footnote{E-mail: $wxm_-gdut$@163.com}\\
$^1\mbox{Department}$ of Mathematics, Guangdong University of Technology,\\
Guangzhou 510520, China\\
$^2 \mbox{School}$ of Mathematics and Statistics,\\
Huazhong University of Science and Technology,\\
Wuhan 430074, China}

\date{}
\maketitle

\begin{abstract}
This paper is devoted to the existence and Lipschitz continuity of global conservative weak solutions in time for the modified two-component Camassa-Holm system on the real line.
We obtain the global weak solutions via a coordinate transformation into the Lagrangian coordinates.
The key ingredients in our analysis are the energy density given by the positive Radon measure
and the proposed new distance functions as well.\\



\end{abstract}

\section{Introduction}
\newtheorem {remark1}{Remark}[section]
\newtheorem{theorem1}{Theorem}[section]
\newtheorem{lemma1}{Lemma}[section]
\newtheorem{definition1}{Definition}[section]

In this paper, we consider the Cauchy problem for the following
modified two-component Camassa-Holm system (M2CH):
\begin{equation}\label{1-M2CH-original}
\left\{\begin{array}{ll}
m_{t}+um_{x}+2mu_x =-\rho\bar{\rho}_{x} ,&t > 0,\,x\in \mathbb{R},\\
\rho_{t}+(\rho u)_x=0, &t > 0,\,x\in \mathbb{R},
\end{array}\right.
\end{equation}
where $m=u-u_{xx}$ and $\rho=(1-\partial^2_{x})(\bar{\rho}-\bar{\rho}_0)$.
System (\ref{1-M2CH-original}) is written in terms of velocity $u$ and locally averaged density $\bar{\rho}$
(or depth, in the shallow-water interpretation) and $\bar{\rho}_0$ is taken to be constant.
As geodesic motion on the semidirect product Lie group with respect to a certain metric,
the system (\ref{1-M2CH-original}) was firstly proposed in \cite{Holm} and proved that
it allows singular solutions in both variables $m$ and $\rho$, not just the fluid momentum.

For $\rho\equiv 0$, system (\ref{1-M2CH-original}) becomes the celebrated Camassa-Holm equation (CH):
\begin{eqnarray*}
m_t+ u m_x+ 2 u_x m=0,\,\ \,\ m=u-u_{xx},
\end{eqnarray*}
which models the unidirectional propagation of shallow water waves over a flat bottom \cite{C-H}.
CH is also a model for the propagation of axially symmetric waves in hyper-elastic rods \cite{Dai}.
It has a bi-Hamiltonian structure and is completely integrable \cite{C-H}.
Its solitary waves are peaked solitons (peakons) \cite{C-H-H,Con-E}, and they are orbitally stable \cite{C-S,C-S1}.
It is noted that the peakons replicate a feature that is characteristic for the waves of great height
-- waves of the largest amplitude that are exact traveling wave solutions of the governing equations for irrotational water waves, cf. \cite{Cinvent,C-Eann}.
The Cauchy problem and initial boundary value problem for CH
have been studied extensively \cite{B-C1,B-C2,C-Ep,C-Ec,E-Y1}.
It has been shown that this equation is locally well-posed \cite{C-Ep,C-Ec}.
Moreover, it has both global strong solutions \cite{Cf,C-Ep,C-Ec} and blow-up solutions within finite time \cite{Cf,C-Ep,C-E,C-Ec}.
It is worthy to point out the advantage of CH in comparison with the KdV equation lies in the fact that
CH has peakons and models wave breaking \cite{C-H-H,C-E}
(namely, the wave remains bounded while its slope becomes unbounded in finite time \cite{Wh}).
Moreover, it possess global weak solutions, see the discussions in \cite{B-C1,B-C2,G-H-Rdcds,H-R1,X-Z}.

The Cauchy problem and initial boundary value problem for system (\ref{1-M2CH-original})
have been investigated in many works, cf. \cite{Guan1,Guan2,guan-yin2,tan-yin1,Yan1,Yan3,Yan2}.
However, in the present paper, we reformulate the considered system to a semilinear system of ODEs
by means of a transformation between Eulerian and Lagrangian coordinates,
which is distinct from those in \cite{Guan2,guan-yin2,tan-yin1}.
As a result, the global existence of conservative weak solution to the system on the real line is constructed.
Moreover, it is noted to point out that we introduce some new distances to derive
the Lipschitz continuity for the obtained weak solution, which, on the other hand, 
implies the uniqueness of the weak solution in some sense.

Now, let us provide the framework in which we shall reformulate system (\ref{1-M2CH-original}).
Set
$\gamma=\bar{\rho}-\bar{\rho}_0$.
By using the identity $(1-\partial^{2}_{x})^{-1}f = p\ast f $
with the Green function $p(x)\triangleq\frac{1}{2}e^{-|x|} (x\in \mathbb{R})$,
one can rewrite the Cauchy problem for system (\ref{1-M2CH-original}) as follows:
\begin{equation}\label{1-M2CH}
\left\{\begin{array}{ll}
u_{t}+uu_{x}=-\partial_x(1-\partial_{x}^{2})^{-1}(u^2+\frac{1}{2}u_x^2+\frac{1}{2}{\gamma}^2-\frac{1}{2}{\gamma_x}^2),
&t> 0,\,x\in \mathbb{R},\\
\gamma_{t}+u\gamma_{x}=-(1-\partial_{x}^{2})^{-1}((u_x\gamma_x)_x+u_x\gamma),
&t> 0,\,x\in \mathbb{R},\\
u(0,x) = u_{0}(x),&x\in \mathbb{R},\\
\gamma(0,x)=\gamma_{0}(x),\ &x\in \mathbb{R}.
\end{array}\right.
\end{equation}

The rest of our paper is organized as follows.
In Section 2, we prove the global existence and Lipschitz continuity of the solutions to system (\ref{1-M2CH})
in Lagrangian coordinates.
In Section 3, we state the stability of the obtained solutions under a new distance
in the setting of Lagrangian coordinates.
In Section 4, the existence of global weak solutions to the system (\ref{1-M2CH}) is proved.
In Section 5, we show the Lipschitz continuity of the weak solution,
which implies the uniqueness of the solution with some certain form.

\bigskip
\noindent \textbf{Notations.}
In the whole paper, we denote by $\ast$ the spatial convolution.
Given a Banach space $Z$, we denote its norm by $\| \cdot\|_{Z}$.
Since all spaces of functions are over $\mathbb{R}$, for simplicity,
we drop $\mathbb{R}$ in our notations of function spaces if there is no ambiguity.

\section{Global and Lipschitz continuous solutions  in Lagrangian coordinates }
\newtheorem {remark2}{Remark}[section]
\newtheorem{theorem2}{Theorem}[section]
\newtheorem{lemma2}{Lemma}[section]
\newtheorem{corollary2}{Corollary}[section]
\newtheorem{definition2}{Definition}[section]

In this section, we discuss the global existence and Lipschitz continuity of the solutions to system (\ref{1-M2CH})
in Lagrangian coordinates. For this,
let us first introduce the spaces $V$ and $V_1$ as follows:
\begin{eqnarray*}
V \triangleq \{f\in C_{b}(\mathbb R)= C(\mathbb R)\cap L^\infty(\mathbb R) \big| f_\xi\in L^2(\mathbb R) \}
\end{eqnarray*}
equipped with the norm $\|f\|_V=\|f\|_{L^\infty}+\|f_\xi\|_{L^2}$, and
\begin{eqnarray*}
V_1 \triangleq \{g\big| g-Id\in V\}.
\end{eqnarray*}

Then the characteristics $y : \mathbb{R}\rightarrow V_1, t\mapsto y(t,\cdot)$ is the solution to
\begin{equation*}\label{2-characteristicsODE}
\left\{\begin{array}{ll}
y_t(t,\xi)=u(t,y(t,\xi)),\\
y(t,\xi)|_{t=0}=y(0,\xi),
\end{array}\right.
\end{equation*}
where $u$ is the first component of the solutions to system (\ref{1-M2CH}).

Set
\begin{align*}
P_1&=p\ast(u^2+\frac{1}{2}u_x^2+\frac{1}{2}{\gamma}^2-\frac{1}{2}{\gamma_x}^2),\\
P_2&=p\ast(u_x\gamma_x),\\
P_3&=p\ast(u_x\gamma).
\end{align*}
For convenience, denote
\begin{equation} \label{2-rewritten-xi}
U(t,\xi)=u(t,y(t,\xi)), \ \ \Gamma(t,\xi)=\gamma(t,y(t,\xi)),\ \  R(t,\xi)=\gamma_x(t,y(t,\xi)),
\end{equation}
\begin{equation*}
P_i(t,\xi)=P_i(t,y(t,\xi)), \ \  Q_i(t,\xi)=P_{i,x}(t,y(t,\xi)),\quad i=1,2,3,
\end{equation*}
and define the Lagrangian energy cumulative distribution as
\begin{equation*} \label{2-Lagrangian-energy}
H(t,\xi)\triangleq\int_{-\infty}^{y(t,\xi)}(u^2+u_x^2+\gamma^2+\gamma_{x}^2)(t,x)dx.
\end{equation*}
Then we can perform the change of variables to write the convolution as an integral
with respect to the new variable $\eta$.
After a straight calculation, one deduces
\begin{align*}
P_1(t,\xi)&=\frac{1}{4}\int_\mathbb R e^{-|y(t,\xi)-y(t,\eta)|}(U^2y_\xi-2R^2y_\xi+H_\xi)d\eta,\\
P_2(t,\xi)&=\frac{1}{2}\int_\mathbb R e^{-|y(t,\xi)-y(t,\eta)|}(RU_\xi)d\eta ,\nonumber\\
P_3(t,\xi)&=\frac{1}{2}\int_\mathbb R e^{-|y(t,\xi)-y(t,\eta)|}(\Gamma U_\xi)d\eta,\nonumber
\end{align*}
and
\begin{align*}
Q_1(t,\xi)&=-\frac{1}{4}\int_\mathbb R sgn((y(t,\xi)-y(t,\eta))  e^{-|y(t,\xi)-y(t,\eta)|}(U^2y_\xi-2R^2y_\xi+H_\xi)d\eta,\\
Q_2(t,\xi)&=-\frac{1}{2}\int_\mathbb R sgn((y(t,\xi)-y(t,\eta)) e^{-|y(t,\xi)-y(t,\eta)|}(RU_\xi)d\eta ,\nonumber\\
Q_3(t,\xi)&=-\frac{1}{2}\int_\mathbb R sgn((y(t,\xi)-y(t,\eta))e^{-|y(t,\xi)-y(t,\eta)|}(\Gamma U_\xi)d\eta.\nonumber
\end{align*}

Observe that $y$ is an increasing function for any fixed $t$, which will be shown later.
Then $P_i$ and $Q_i$ have the following equivalent forms:
\begin{align} \label{2-P-equivalent-forms}
P_1(t,\xi)&=\frac{1}{4}\int_\mathbb R e^{-sgn(\xi-\eta)(y(t,\xi)-y(t,\eta))}(U^2y_\xi-2R^2y_\xi+H_\xi)d\eta,\\
P_2(t,\xi)&=\frac{1}{2}\int_\mathbb R e^{-sgn(\xi-\eta)(y(t,\xi)-y(t,\eta))}(RU_\xi)d\eta ,\nonumber\\
P_3(t,\xi)&=\frac{1}{2}\int_\mathbb R e^{-sgn(\xi-\eta)(y(t,\xi)-y(t,\eta))}(\Gamma U_\xi)d\eta,\nonumber
\end{align}
and
\begin{align} \label{2-Q-equivalent-forms}
Q_1(t,\xi)&=-\frac{1}{4}\int_\mathbb R sgn(\xi-\eta)  e^{-sgn(\xi-\eta)(y(t,\xi)-y(t,\eta))}(U^2y_\xi-2R^2y_\xi+H_\xi)d\eta,\\
Q_2(t,\xi)&=-\frac{1}{2}\int_\mathbb R sgn(\xi-\eta) e^{-sgn(\xi-\eta)(y(t,\xi)-y(t,\eta))}(RU_\xi)d\eta ,\nonumber\\
Q_3(t,\xi)&=-\frac{1}{2}\int_\mathbb R sgn(\xi-\eta)e^{-sgn(\xi-\eta)(y(t,\xi)-y(t,\eta))}(\Gamma U_\xi)d\eta.\nonumber
\end{align}
Besides, we can check that the following equalities hold:
\begin{align} \label{2-PQ-equalities}
P_{1,\xi}=Q_1y_\xi,\ \ \ \ &Q_{1,\xi}=-\frac{1}{2}(U^2y_\xi-2R^2y_\xi+H_\xi)+P_1y_\xi,\\
P_{2,\xi}=Q_2y_\xi,\ \ \ \ &Q_{2,\xi}=-RU_\xi+P_2y_\xi ,\nonumber\\
P_{3,\xi}=Q_3y_\xi,\ \ \ \ &Q_{3,\xi}=-\Gamma U_\xi+P_3y_\xi.\nonumber
\end{align}

Next, we introduce another new variable $\zeta(t,\xi)=y(t,\xi)-\xi$,
then a new system based on system (\ref{1-M2CH}) is derived as follows:
\begin{equation} \label{2-new-system}
\left\{\begin{array}{ll}
\zeta_{t}=U,\\
U_t=-Q_1,\\
\Gamma_t=-Q_2-P_3,\\
R_t =-P_2-Q_3,\\
H_t=U^3-2P_1U-2P_2\Gamma-2Q_3\Gamma.
\end{array}\right.
\end{equation}
Define
$$ E \triangleq V\times H^{1}\times H^{1}\times(L^\infty\cap L^2)\times V$$
with the norm $\|X\|_E=\|\zeta\|_V+\|U\|_{H^1}+\|\Gamma\|_{H^1}+\|R\|_{L^\infty}+\|R\|_{L^2}+\|H\|_V.$
Then we have the following Lipschitz estimates for $P_i$ and $Q_i$ $(i=1,2,3)$.
\begin{lemma2}
For any $X=(\zeta,U,\Gamma,R,H)\in E$, define the maps $F_i(X)=P_i$ and $G_i(X)=Q_i (i=1,2,3)$,
where $P_i$ and $Q_i$ are given by (\ref{2-P-equivalent-forms}) and (\ref{2-Q-equivalent-forms}) with $y=\zeta+Id$, respectively.
Then $F_i$ and $G_i$ are $B$-Lipschitz from $E$ to $H^1$,
namely, they are Lipschitz continuous from the bounded sets in $E$ to $H^1$.
More precisely, let $B_M$ be the closed ball with radius $M$ in $E$.
Then for any $X, \tilde{X}\in E$, we have
\begin{align} \label{2-Lipschitz-estimates-F}
\|F_i(X)-F_i(\tilde{X})\|_{H^1}\leq C_M\|X-\tilde{X}\|_E
\end{align}
and
\begin{align*}
\|G_i(X)-G_i(\tilde{X})\|_{H^1}\leq C_M\|X-\tilde{X}\|_E,
\end{align*}
where $i=1,2,3$ and the positive constant $C_M$ depends only on $M$.
\end{lemma2}

\begin{proof}
We here only prove the estimate for $F_1$,  the others can be handled in a similar way.
Indeed, from (\ref{2-P-equivalent-forms}),  we rewrite $F_1$ as
\begin{align*}
F_1(X)(\xi)=&\frac{e^{-\zeta(\xi)}}{4}\int_\mathbb R\chi\{\eta<\xi\}e^{-|\xi-\eta|}e^{\zeta(\eta)}[(U^2-2R^2)(1+\zeta_\xi)+H_\xi](\eta)d\eta\\
&+\frac{e^{\zeta(\xi)}}{4}\int_\mathbb R\chi\{\eta>\xi\}e^{-|\xi-\eta|}e^{-\zeta(\eta)}[(U^2-2R^2)(1+\zeta_\xi)+H_\xi](\eta)d\eta\\
\triangleq &I_1+I_2,
\end{align*}
where $\chi_A$ is the indicator function of some set $A$.

Let $h(\xi)=\chi_{\{\xi>0\}}(\xi)e^{-\xi}$ and define the map $P:v\mapsto h\ast v.$
Thanks to Lemma 2.1 in \cite{H-R1}, the map $P$ is continuous from $L^2$ into $H^1.$
If we denote
$$R(X)\triangleq e^\zeta((U^2-2R^2)(1+\zeta_\xi)+H_\xi),$$
then $I_1=\frac{e^{-\zeta(\xi)}}{4}P\circ R(X)(\xi).$

Now we check that $R$ is $B$-Lipschitz from $E$ to $L^2.$ Indeed,
for any $X, \tilde{X}\in B_M$,
by using $|e^x-e^y|\leq e^{max\{|x|,|y|\}}|x-y|$
and $\|f\|_{L^\infty}\leq\|f\|_{H^1}$, one infers
\begin{align*}
&\|R(X)-R(\tilde{X})\|_{L^2}\\
=&\|e^\zeta [(U^2-2R^2)(1+\zeta_\xi)+H_\xi]-e^{\tilde{\zeta}}[(\tilde{U}^2-2\tilde{R}^2)(1+\tilde{\zeta}_\xi)+\tilde{H}_\xi]\|_{L^2}\\
\leq&\|e^\zeta U^2-e^{\tilde{\zeta}}\tilde{U}^2\|_{L^2}+2\|e^\zeta R^2-e^{\tilde{\zeta}}\tilde{R}^2\|_{L^2}+\|e^\zeta U^2\zeta_\xi-e^{\tilde{\zeta}}\tilde{U}^2\tilde{\zeta}_\xi\|_{L^2}\\
&+2\|e^\zeta R^2\zeta_\xi-e^{\tilde{\zeta}}\tilde{R}^2\tilde{\zeta}_\xi\|_{L^2}+\|e^\zeta H_\xi-e^{\tilde{\zeta}}\tilde{H}_\xi\|_{L^2}\\
\leq&\|e^\zeta U^2-e^{\zeta}\tilde{U}^2\|_{L^2}+\|e^{\zeta}\tilde{U}^2-e^{\tilde{\zeta}}\tilde{U}^2\|_{L^2}+2\|e^\zeta R^2-e^{\zeta}\tilde{R}^2\|_{L^2 }\\
&+2\|e^{\zeta}\tilde{R}^2-e^{\tilde{\zeta}}\tilde{R}^2\|_{L^2}+\|e^\zeta U^2\zeta_\xi-e^{\zeta}U^2\tilde{\zeta}_\xi\|_{L^2}
+\|e^{\zeta}U^2\tilde{\zeta}_\xi-e^{\zeta}\tilde{U}^2\tilde{\zeta}_\xi\|_{L^2}\\
&+\|e^{\zeta}\tilde{U}^2\tilde{\zeta}_\xi-e^{\tilde{\zeta}}\tilde{U}^2\tilde{\zeta}_\xi\|_{L^2}
+2\|e^\zeta R^2\zeta_\xi-e^{\zeta}R^2\tilde{\zeta}_\xi\|_{L^2}
+2\|e^{\zeta}R^2\tilde{\zeta}_\xi-e^{\zeta}\tilde{R}^2\tilde{\zeta}_\xi\|_{L^2}\\
&+2\|e^{\zeta}\tilde{R}^2\tilde{\zeta}_\xi-e^{\tilde{\zeta}}\tilde{R}^2\tilde{\zeta}_\xi\|_{L^2}
+\|e^{\zeta}H_\xi-e^\zeta\tilde{H}_\xi\|_{L^2}+\|e^\zeta\tilde{H}_\xi-e^{\tilde{\zeta}}\tilde{H}_\xi\|_{L^2}\\
\leq&\|e^\zeta\|_{L^\infty}\|U+\tilde{U}\|_{L^\infty}\|U-\tilde{U}\|_{L^2}+\|e^\zeta-e^{\tilde{\zeta}}\|_{L^\infty}\|\tilde{U}\|_{L^\infty}\|\tilde{U}\|_{L^2}\\
&+2\|e^\zeta\|_{L^\infty}\|R+\tilde{R}\|_{L^\infty}\|R-\tilde{R}\|_{L^2}+2\|e^\zeta-e^{\tilde{\zeta}}\|_{L^\infty}\|\tilde{R}\|_{L^\infty}\|\tilde{R}\|_{L^2}\\
&+\|e^\zeta\|_{L^\infty}\|U\|^2_{L^\infty}\|\zeta_\xi-\tilde{\zeta}_\xi\|_{L^2}
+\|e^\zeta\|_{L^\infty}\|\tilde{\zeta}_\xi\|_{L^2}\|U+\tilde{U}\|_{L^\infty}\|U-\tilde{U}\|_{L^2}\\
&+\|e^\zeta-e^{\tilde{\zeta}}\|_{L^\infty}\|\tilde{U}\|^2_{L^\infty}\|\tilde{\zeta}_\xi\|_{L^2}
+2\|e^\zeta\|_{L^\infty}\|R\|^2_{L^\infty}\|\zeta_\xi-\tilde{\zeta}_\xi\|_{L^2}\\
&+2\|e^\zeta\|_{L^\infty}\|\tilde{\zeta}_\xi\|_{L^2}\|R+\tilde{R}\|_{L^\infty}\|R-\tilde{R}\|_{L^2}
+2\|e^\zeta-e^{\tilde{\zeta}}\|_{L^\infty}\|\tilde{R}\|^2_{L^\infty}\|\tilde{\zeta}_\xi\|_{L^2}\\
&+\|e^\zeta\|_{L^\infty}\|H_\xi-\tilde{H}_\xi\|_{L^2}+\|e^\zeta-e^{\tilde{\zeta}}\|_{L^\infty}\|\tilde{H}_\xi\|_{L^2}\\
\leq&2Me^M\|U-\tilde{U}\|_{L^2}+M^2e^M\|\zeta-\tilde{\zeta}\|_{L^\infty}+4Me^M\|R-\tilde{R}\|_{L^2}\\
&+2M^2e^M\|\zeta-\tilde{\zeta}\|_{L^\infty}+M^2e^M\|\zeta_\xi-\tilde{\zeta}_\xi\|_{L^2}+2M^2e^M\|U-\tilde{U}\|_{L^2}\\
&+M^3e^M\|\zeta-\tilde{\zeta}\|_{L^\infty}+2M^2e^M\|\zeta_\xi-\tilde{\zeta}_\xi\|_{L^2}+4M^2e^M\|R-\tilde{R}\|_{L^2}\\
&+2M^3e^M\|\zeta-\tilde{\zeta}\|_{L^\infty}+e^M\|H_\xi-\tilde{H}_\xi\|_{L^2}+Me^M\|\zeta-\tilde{\zeta}\|_{L^\infty}\\
\leq&C(M)\|X-\tilde{X}\|_E,
\end{align*}
where $C(M)=(1+7M+12M^2+3M^3)e^M.$

For any bounded set $B$ in $E$, there is an $M>0$ such that $B\subset B_M$.
Thus, $R$ is $B$-Lipschitz from $E$ into $L^2$.
Note that $P$ is a continuous linear operator from $L^2$ to $H^1$.
Then $I_1=P\circ P$ is $B$-Lipschitz from $E$ to $H^1.$
Similarly, so is $I_2$.
Hence, $F_1$ is a Lipschitz map from $E$ to $H^1$ and satisfies (\ref{2-Lipschitz-estimates-F}).
Therefore, we have proven the lemma.
\end{proof}

\begin{lemma2}
For any initial data $\bar{X}\in E$, there exists a time $T=T(\|\bar{X}\|_E)>0$ such that
the system (\ref{2-new-system}) admits a unique solution in $C^1([0,T],E)$.
\end{lemma2}

\begin{proof}
Define the map $F: E\rightarrow E$ by
$$F(X)=(U,-Q_1,-Q_2-P_3,-P_2-Q_3,U^3-2P_1U-2P_2\Gamma-2Q_3\Gamma).$$
Then Lemma 2.1 and the standard ODE theory on Banach spaces yield the desired result.
\end{proof}

Now, we turn our attention to the global existence of the solutions to system (\ref{2-new-system}).
Here we are interested in a special class of initial data which belong to
$E_0 \triangleq (W^{1,\infty}\times W^{1,\infty}\times W^{1,\infty}\times L^\infty\times W^{1,\infty})\cap E$.

For any $X_0\in E_0$, thanks to Lemma 2.2,
the system (\ref{2-new-system}) has a unique short time solution $X\in C([0,T],E)$.
Let $U,\Gamma, P_i,Q_i\in C([0,T],H^1)$ and $R\in C([0,T],L^\infty \cap L^2)$.
For any fixed $\xi\in\mathbb R$,  we can solve the following system
\begin{equation} \label{2-new-system-2}
\left\{\begin{array}{ll}
\alpha_t=\beta,\\
 \beta_t=\frac{1}{2}(U^2-2R^2-2p_1)(1+\alpha)+\frac{1}{2}\delta,\\
\kappa_t =-(P_2+Q_3)(1+\alpha)+R\beta,\\
\delta_t=(3U^2-2P_1+2\Gamma^2)\beta-(2UQ_1+2Q_2\Gamma+2P_3\Gamma)(1+\alpha)-(2P_2+2Q_3)\kappa, \end{array}\right.
\end{equation}
by substituting $\zeta_\xi,U_\xi,\Gamma_\xi$ and $H_\xi$ in system (\ref{2-new-system})
by $\alpha,\beta,\kappa$ and $\delta$, respectively.

Similar to the proof of Lemma 3.4 in \cite{tan-yin1}, one can readily get the following result.
\begin{lemma2}
Let $X\in C([0,T],E)$ be the solution to system (\ref{2-new-system}) with the initial data $\bar{X}\in E_0$.
Then $(\alpha,\beta,\kappa,\delta)$ solves system (\ref{2-new-system-2}) for any fixed $\xi.$
Moreover, for all $t\in[0,T]$ and almost every $\xi\in\mathbb R$, we have
\begin{align} \label{2-equality-two-systems}
(\alpha(t,\xi),\beta(t,\xi),\kappa(t,\xi),\delta(t,\xi))=(\zeta_\xi(t,\xi),U_\xi(t,\xi),\Gamma_\xi(t,\xi),H_\xi(t,\xi)).
\end{align}
\end{lemma2}

\begin{definition2}
The set $\mathbb G$ is composed of all $(\zeta, U,  \Gamma, R, H)\in E$,  such that
\begin{align}
&(\zeta, U,\Gamma, H)\in( W^{1,\infty})^4, \label{2-condition-1} \\
&y_\xi\geq0, \ \ H_\xi\geq0, \ y_\xi+H_\xi>0, \ \textrm{a.e.\ on} \ \mathbb{R}, \,\  \textrm{and} \
\lim_{\xi\rightarrow -\infty}H(\xi)=0, \label{2-condition-2}\\
&Ry_\xi=\Gamma_\xi,\ \textrm{a.e.\ on} \ \mathbb{R}, \label{2-condition-3}\\
&y_\xi H_\xi=y^2_\xi (U^2+\Gamma^2+R^2)+U^2_\xi, \ \textrm{a.e.\ on} \ \mathbb{R}, \,\
\textrm{where}\,\ y(\xi)=\zeta(\xi)+\xi. \label{2-condition-4}
\end{align}
\end{definition2}

According to Lemma 2.3, one can prove the following useful lemma.
\begin{lemma2}
The set $\mathbb G$ is preserved by system (\ref{2-new-system}). That is,
if the initial data $X_0\in\mathbb G$, then the corresponding solution $X(t,\cdot)$ to system (\ref{2-new-system})
also belongs to $\mathbb G$ for all $t\in[0,T]$.
Furthermore, $y_\xi(t,\xi)>0$ holds true for almost every $(t,\xi)\in[0,T]\times \mathbb{R}$.
\end{lemma2}

\begin{proof}
By Lemma 2.3, (\ref{2-condition-1}) holds for all $t\in[0,T]$.
Next, for any fixed $\xi$ which satisfies $|X(\xi)|\leq\|X\|_{L^\infty}$,
by system (\ref{2-new-system}), we have $(Ry_\xi)_t=\Gamma_{\xi t}$.
This leads to (\ref{2-condition-3}), provided that $R(0)y_\xi(0)=\Gamma_\xi(0)$.

From the system (\ref{2-new-system-2}) and (\ref{2-equality-two-systems}), one has
\begin{eqnarray*}
(y_\xi H_\xi)_t
&=&y_{\xi t}H_\xi+H_{\xi t}y_\xi\\
&=&U_\xi H_\xi+y_\xi((3U^2-2P_1+2\Gamma^2)U_\xi
-(2UQ_1+2Q_2\Gamma+2P_3\Gamma)y_\xi\\
&&-(2P_2+2Q_3)\Gamma_\xi),
\end{eqnarray*}
and
\begin{align*}
&(y^2_\xi (U^2+\Gamma^2+R^2)+U^2_\xi)_t\\
=&2(UU_t+\Gamma\Gamma_t+RR_t)y_{\xi}^2+2(U^2+\Gamma^2+R^2)y_\xi y_{\xi t}+2U_\xi U_{\xi t}\\
=&-2UQ_1y_\xi^2+2U^2y_\xi U_\xi-2\Gamma(Q_2+P_3)y_\xi^2+2\Gamma^2y_\xi U_\xi-2R(P_2+Q_3)y^2_\xi\\
&+2R^2y_\xi U_\xi+2\{\frac{1}{2}(U^2y_\xi-2R^2y_\xi+H_\xi)-P_1y_\xi\}U_\xi\\
=&U_\xi H_\xi+3U^2U_\xi y_\xi-2P_1y_\xi U_\xi+2\Gamma^2y_\xi U_\xi-2UQ_1y_\xi^2-2\Gamma Q_2y_\xi^2\\
&-2\Gamma P_3y_\xi^2-2(P_2+Q_3)Ry_\xi^2,
\end{align*}
which together with (\ref{2-condition-3}) implies (\ref{2-condition-4}).

Finally, we prove (\ref{2-condition-2}). For this, define
$$ t^*\triangleq \sup\{t\in[0,T] \big| y_\xi(t')\geq0  \,\ \textrm{for\ all}\,\ t'\in[0,t]\}.$$
We claim that $t^*=T$. Suppose not, i.e. $t^*<T$.
Since $y_\xi(t)$ is continuous with respect to $t$, it follows that
\begin{align} \label{2-lemma4-1}
y_\xi(t^*)=0,
\end{align}
which along with (\ref{2-condition-4}) gives $U_\xi(t^*)=0$.
Thanks to system (\ref{2-new-system-2}) again, one deduces
\begin{align}\label{2-lemma4-2}
y_{\xi t}(t^*)=U_\xi(t^*)=0.
\end{align}
By system (\ref{2-new-system}), together with (\ref{2-lemma4-1}) and (\ref{2-lemma4-2}), we get
\begin{align}\label{2-lemma4-3}
y_{\xi tt}(t^*)=U_{\xi t}(t^*)=\frac{1}{2}H_\xi(t^*).
\end{align}
Now, we can show $H_\xi(t^*)>0$. As a matter of fact,
if $H_\xi(t^*)=0,$ then (\ref{2-condition-3}) and (\ref{2-condition-4}) ensure
$$y_\xi(t^*)=U_\xi(t^*)=\Gamma_\xi(t^*)=H_\xi(t^*)=0,$$
which is a contradiction to the uniqueness of the solution to system (\ref{2-new-system-2}).
If $H_\xi(t^*)<0,$ then according to (\ref{2-lemma4-3}), $y_{\xi tt}<0$.
So, $y_\xi(t^*)$ is the strict maximum of $y_\xi$, which also contradicts the definition of $t^*.$
Hence, $H_\xi(t^*)>0$. Then $y_\xi(t^*)$ is the minimum of $y_\xi$. 
Thus, from the definition of $t^*$ again, we have $t^*=T.$
This is a contradiction to the assumption $t^*<T$.
Therefore, $y_\xi(t)\geq 0$ holds for all $t\in[0,T]$.

Furthermore, we have $y_\xi(t)> 0$ for all $t\in[0,T]$. Indeed,
if there is some $t'\in [0,t]$ such that $y_\xi(t')=0$, then $H_\xi(t')<0$ as above arguments.
So, $y_\xi(t')=0$ is the strict minimum of $y_\xi$,
which contradicts the fact $y_\xi(t)\geq 0$ for all $t\in[0,T]$.
Note that (\ref{2-condition-4}) implies $H_\xi(t)\geq 0$ whenever $y_\xi(t)>0$.
Thus, $y_\xi(t)+H_\xi(t)>0$ for $t\in [0,T].$
Therefore, we complete the proof of the lemma.
\end{proof}

With Lemmas 2.2-2.4 in hand, we conclude this section with the following main theorem.
%In what follows, we always denote $y(t,\xi)=\zeta(t,\xi)+\xi.$

\begin{theorem2}
For any initial data $X_0=(\zeta_0, U_0, \Gamma_0, R_0,H_0)\in\mathbb G$,
the system (\ref{2-new-system}) has a unique global solution
$X(t)=(\zeta(t), U(t), \Gamma(t), R(t),H(t))\in\mathbb G$ for all $t>0$, and $X(t)\in C^1(\mathbb R_+;E)$.
Moreover, define the mapping $S_t: \mathbb G\times\mathbb R_+\rightarrow\mathbb G$ as  $S_t(X)=X(t).$
Then the mapping $S_t$ is a continuous semigroup.
Furthermore, let $M,\,T>0$ and set
$$B_M \triangleq \{X=(\zeta, U,\Gamma, R,H)\in\mathbb G \big| X\in E, \|X\|_E\leq M\}.$$
Then there exists a $C=C(M,T)>0$ such that
\begin{equation} \label{2-thm-estimate}
\|S_t(X_\alpha)-S_t(X_\beta)\|_E\leq C \|X_\alpha-X_\beta\|_E,
\end{equation}
for any $X_\alpha, X_\beta\in B_M$.
\end{theorem2}

\begin{proof}
By Lemmas 2.3-2.4 and a contraction argument,
we obtain a short time solution $X(t)\in\mathbb G$ to system (\ref{2-new-system})
with initial data $X_0\in\mathbb G$. In addition,
the solution has a finite maximal existence time $T$ if and only if
$$\lim_{t\rightarrow T}\|X(t)\|_E=\infty.$$

For any $T_1<T$ and $t\in [0,T_1]$, in view of Lemma 2.4 and (\ref{2-condition-2}),
$H(t,\xi)$ is an increasing function with respect to $\xi$.
Hence, the limits $H(t,\pm\infty)\triangleq \lim\limits_{\xi\rightarrow\pm\infty} H(t,\xi)$ exist.
Note that $U(t,\cdot),\, \Gamma(t,\cdot)\in H^1$. Then we have
\begin{equation} \label{2-thm-limit}
\lim_{\xi\rightarrow\pm\infty}U(t,\xi)=\lim_{\xi\rightarrow\pm\infty}\Gamma(t,\xi)=0.
\end{equation}
So, the last equation in system (\ref{2-new-system}) leads to
\begin{align} \label{2-thm-H}
H(t,\xi)=H_0(\xi)+\int_0^t(U^3-2P_1U-2P_2\Gamma-2Q_3\Gamma)(s,\xi)ds.
\end{align}
Since Lemma 2.1 implies that $U,\Gamma,P_1,P_2$ and $Q_3$ are bounded in $L^\infty([0,T]\times\mathbb R)$,
it follows from (\ref{2-thm-limit}), (\ref{2-thm-H}) and the Lebesgue dominated convergence theorem that $H(t,\pm\infty)=H_0(\pm\infty)$ holds for all $t\in[0,T_1].$
Observe that both $H(t,\xi)$ and $H_0(\xi)$ are increasing with respect to $\xi$.
Thus, $\sup_{t\in[0,T)}\|H(t,\cdot)\|_{L^\infty}$ can be bounded by $\|H_0\|_{L^\infty}$.

On the other hand, by using Lemma 2.4 and (\ref{2-condition-4}), one has
\begin{align*}
U^2(t,\xi)=&2\int_{-\infty}^\xi U(t,\eta)U_\xi(\eta)d\eta\\
=&\nonumber2\int_{\eta\leq\xi|y_\xi(\eta)>0}U(t,\eta)U_\xi(\eta)d\eta\\
\leq&\nonumber2\int_{\eta\leq\xi|y_\xi(\eta)>0}|U(t,\eta)U_\xi(\eta)|d\eta\\
\leq&\nonumber\int_{\eta\leq\xi|y_\xi(\eta)>0}\left(U^2(t,\eta)y_\xi(t,\eta)+\frac{U_\xi^2(t,\eta)}{y_\xi(t,\eta)}\right)d\eta\\
\leq&\nonumber\int_{\eta\leq\xi|y_\xi(\eta)>0}H_\xi(\eta)d\eta\\
\leq& H(t,\xi),
\end{align*}
which implies that  $\sup_{t\in[0,T)}\|U(t,\cdot)\|_{L^\infty}$ can be bounded by $\sqrt{\|H_0\|_{L^\infty}}$.\\
Thanks to Lemma 2.4 and (\ref{2-condition-4}) again, we have
\begin{align*}
U^2y_\xi+\Gamma^2y_\xi+R^2y_\xi+\frac{U^2_\xi}{y_\xi}=H_\xi.
\end{align*}
By recalling (\ref{2-P-equivalent-forms}) and (\ref{2-Q-equivalent-forms}), one gets
\begin{align*}
\sup_{t\in[0,T)}\|P_1(t,\cdot)\|_{L^\infty}\leq\sup_{t\in[0,T)}\int_\mathbb R H_\xi(t,\eta)d\eta\leq\|H_0\|_{L^\infty},
\end{align*}
and
\begin{align*}
\sup_{t\in[0,T)}\|Q_1(t,\cdot)\|_{L^\infty}\leq\sup_{t\in[0,T)}\int_\mathbb R H_\xi(t,\eta)d\eta\leq\|H_0\|_{L^\infty}.
\end{align*}
Note that
$|RU_\xi|\leq\frac{1}{2}(R^2y_\xi+\frac{U^2_\xi}{y_\xi})\leq\frac{1}{2} H_\xi$ and
$|\Gamma U_\xi|\leq\frac{1}{2}(\Gamma^2y_\xi+\frac{U^2_\xi}{y_\xi})\leq\frac{1}{2} H_\xi$.
We deduce
\begin{align*}
\sup_{t\in[0,T)}\|P_2(t,\cdot)\|_{L^\infty}\leq\sup_{t\in[0,T)}\int_\mathbb R H_\xi(t,\eta)d\eta\leq\|H_0\|_{L^\infty},
\end{align*}
\begin{align*}
\sup_{t\in[0,T)}\|Q_2(t,\cdot)\|_{L^\infty}\leq\sup_{t\in[0,T)}\int_\mathbb R H_\xi(t,\eta)d\eta\leq\|H_0\|_{L^\infty},
\end{align*}
\begin{align*}
\sup_{t\in[0,T)}\|P_3(t,\cdot)\|_{L^\infty}\leq\sup_{t\in[0,T)}\int_\mathbb R H_\xi(t,\eta)d\eta\leq\|H_0\|_{L^\infty},
\end{align*}
and
\begin{align*}
\sup_{t\in[0,T)}\|Q_3(t,\cdot)\|_{L^\infty}\leq\sup_{t\in[0,T)}\int_\mathbb R H_\xi(t,\eta)d\eta\leq\|H_0\|_{L^\infty}.
\end{align*}
From system (\ref{2-new-system}),  we obtain
\begin{align*}
\sup_{t\in[0,T)}\|\zeta(t)\|_{L^\infty}\leq\|\zeta_0\|_{L^\infty}+\int_0^t\|U(s)\|_{L^\infty}ds
\leq\|\zeta_0\|_{L^\infty}+T\sqrt{\|H_0\|_{L^\infty}},
\end{align*}
\begin{align*}
\sup_{t\in[0,T)}\|\Gamma(t)\|_{L^\infty}\leq\|\Gamma_0\|_{L^\infty}+\int_0^t\|(Q_2+P_3)(s)\|_{L^\infty}ds
\leq\|\Gamma_0\|_{L^\infty}+2T\|H_0\|_{L^\infty},
\end{align*}
and
\begin{align*}
\sup_{t\in[0,T)}\|R(t)\|_{L^\infty}\leq\|R_0\|_{L^\infty}+\int_0^t\|(P_2+Q_3)(s)\|_{L^\infty}ds
\leq\|R_0\|_{L^\infty}+2T\|H_0\|_{L^\infty}.
\end{align*}

Denote
\begin{eqnarray*}
C_1
&\triangleq&\sup_{t\in[0,T)}\{\|\zeta(t,\cdot)\|_{L^\infty}+\|U(t,\cdot)\|_{L^\infty}
+\|\Gamma(t,\cdot)\|_{L^\infty}+\|R(t,\cdot)\|_{L^\infty}\\
&&+\sum\limits_{i=1}^3(\|P_i(t,\cdot)\|_{L^\infty}+\|Q_i(t,\cdot)\|_{L^\infty})\}.
\end{eqnarray*}
Then $C_1$ is finite and depends only on $T$ and the initial data.
Similar to the proof of Lemma 2.1, for $i=1,2,3$, we have
\begin{align} \label{2-thm-estimate-PQ}
&\|P_i(t\cdot)\|_{L^2},\|Q_i(t,\cdot)\|_{L^2}\\
\leq &C(\|U\|_{L^2}+\|\Gamma\|_{L^2}+\|R\|_{L^2}+\|\zeta_\xi\|_{L^2}+\|U_\xi\|_{L^2}+\|H_\xi\|_{L^2}),\nonumber
\end{align}
where $C$ is also finite and depends only on $C_1$.
Then from the system (\ref{2-new-system}) and (\ref{2-thm-estimate-PQ}), one gets
\begin{align*}
&\|R(t,\cdot)\|_{L^2}\\
\leq &\|R(t,\cdot)\|_{L^2}+C\int_0^t(\|U\|_{L^2}+\|\Gamma\|_{L^2}+\|R\|_{L^2}
+\|\zeta_\xi\|_{L^2}+\|U_\xi\|_{L^2}+\|H_\xi\|_{L^2})(s)ds,\nonumber
\end{align*}
which together with Lemma 2.3 and system (\ref{2-new-system-2}) lead to
\begin{eqnarray*}
Z(t)\leq Z(0)+C\int_0^tZ(\tau)d\tau,
\end{eqnarray*}
where $Z(t)\triangleq \|\zeta_\xi\|_{L^2}+\|U_\xi\|_{L^2}+\|\Gamma_\xi\|_{L^2}+\|H_\xi\|_{L^2}+\|R\|_{L^2}$
and the positive constant $C$ depends only on $C_1$.

Taking advantage of the Gronwall inequality, one deduces that $\sup_{t\in[0,T)}\|X(t)\|_E$ is finite.
Then the standard ODE theory implies that $S_t$ is a continuous semigroup.
And thus, we have obtained the global existence of the solutions.

In order to prove the theorem, it suffices to show (\ref{2-thm-estimate}). Indeed,
for any $X_\alpha, X_\beta\in B_M,$ we see $X_\alpha- X_\beta,\, S_t(X_\alpha)-S_t(X_\beta)\in E$
and
\begin{align*}
\|S_t(X)\|_E\leq C(T,\|X\|_E),
\end{align*}
for all $t\in[0,T]$ with any $T>0$.

By the second equation in system (\ref{2-new-system}), we get
\begin{align*}
U_\alpha(t,\xi)-U_\beta(t,\xi)=U_\alpha(\xi)-U_\beta(\xi)+\int_0^t(Q_\alpha-Q_\beta)(s,\xi)ds,
\end{align*}
which yields
\begin{align*}
\|U_\alpha(t,\cdot)-U_\beta(t,\cdot)\|_{L^2}
\leq\|U_\alpha-U_\beta\|_{L^2}+\int_0^t\|Q_\alpha(s)-Q_\beta(s)\|_{L^2}ds,
\end{align*}
and
\begin{align*}
\|(U_\alpha(t,\cdot)-U_\beta(t,\cdot))_\xi\|_{L^2}
\leq\|(U_\alpha-U_\beta)_\xi\|_{L^2}+\int_0^t\|(Q_\alpha(s)-Q_\beta(s))_\xi\|_{L^2}ds.
\end{align*}
Then the above two inequalities and Lemma 2.1 imply
\begin{align} \label{2-thm-U-H1}
\|U_\alpha(t,\cdot)-U_\beta(t,\cdot)\|_{H^1}
\leq\|U_\alpha-U_\beta\|_{H^1}+C\int_0^t\|X_\alpha(s)-X_\beta(s)\|_{E}ds,
\end{align}
where $C=C(T,M)>0$.\\
Likewise, from the equations $(\ref{2-new-system})_3$ and $(\ref{2-new-system})_4$, we have
\begin{align} \label{2-thm-Gamma-H1}
\|\Gamma_\alpha(t,\cdot)-\Gamma_\beta(t,\cdot)\|_{H^1}
\leq\|\Gamma_\alpha-\Gamma_\beta\|_{H^1}+C\int_0^t\|X_\alpha(s)-X_\beta(s)\|_{E}ds,
\end{align}
and
\begin{align} \label{2-thm-R-L-L2}
\|R_\alpha(t,\cdot)-R_\beta(t,\cdot)\|_{L^\infty \cap L^2}
\leq \|R_\alpha-R_\beta\|_{L^\infty\cap L^2}+C\int_0^t\|X_\alpha(s)-X_\beta(s)\|_{E}d s.
\end{align}
By virtue of $(\ref{2-new-system})_1$, $(\ref{2-new-system})_3$ and (\ref{2-thm-U-H1}), we get
\begin{align} \label{2-thm-zta-V}
\|\zeta_\alpha(t,\cdot)-\zeta_\beta(t,\cdot)\|_{V}
\leq \|\zeta_\alpha-\zeta_\beta\|_{V}+C\int_0^t\|X_\alpha(s)-X_\beta(s)\|_{E}d s.
\end{align}
From the last equation in system (\ref{2-new-system}), we infer
\begin{eqnarray*}
H_\alpha(t,\xi)-H_\beta(t,\xi)
&=& H_\alpha(\xi)-H_\beta(\xi)+\int_0^t((U^3-2P_1U-2P_2\Gamma-2Q_3\Gamma)_\alpha\\
&&-(U^3-2P_1U-2P_2\Gamma-2Q_3\Gamma)_\beta)(s,\xi)d s,
\end{eqnarray*}
which together with (\ref{2-thm-U-H1})-(\ref{2-thm-zta-V}) and Lemma 2.1 yield
\begin{align*}
\|H_\alpha(t,\cdot)-H_\beta(t,\cdot)\|_{V}
\leq \|H_\alpha-H_\beta\|_{V}+C_M\int_0^t\|X_\alpha(s)-X_\beta(s)\|_E d s.
\end{align*}

Hence, we obtain
\begin{align*}
\|X_\alpha(t,\cdot)-X_\beta(t,\cdot)\|_{E}
\leq\|X_\alpha-X_\beta\|_{E}+ C_M\int_0^t\|X_\alpha(s)-X_\beta(s)\|_{E}ds,
\end{align*}
where $C=C(M,T)>0$. Making use of the Gronwall inequality again, one reaches (\ref{2-thm-estimate}).
Therefore, we complete the proof of the theorem.
\end{proof}

\section{Stability of the solutions under a new distance}
\newtheorem {remark3}{Remark}[section]
\newtheorem{theorem3}{Theorem}[section]
\newtheorem{lemma3}{Lemma}[section]
\newtheorem{definition3}{Definition}[section]

In this section, we investigate the stability of the weak solutions to system (\ref{1-M2CH}) under a new distance.
To this end, let us first denote by G the subgroup of the group of homeomorphisms
from $\mathbb R$ to $\mathbb R$ as follows:
$$G \triangleq \{f \big| f-Id\in W^{1,\infty}, f^{-1}-Id\in W^{1,\infty}, f_\xi-1\in L^2\}.$$
For any $\alpha>0,$ we introduce the subsets $G_\alpha$ of $G$ by
$$G_\alpha \triangleq \{f\in G \big| \|f-Id\|_{W^{1,\infty}}+\|f^{-1}-Id\|_{W^{1,\infty}}\leq\alpha\}.$$

\begin{lemma3}
\cite{H-R1}
If $f\in G_\alpha,$ then $\frac{1}{1+\alpha}\leq f_\xi\leq 1+\alpha$ almost everywhere.
Conversely, if f is absolutely continuous, $f-Id\in L^\infty$
and there exists some $c\geq1$ such that $\frac{1}{c}\leq f_\xi\leq c$ almost everywhere,
then $f\in G_\alpha$ for some $\alpha$ depending only on $c$ and $\|f-Id\|_{L^\infty}.$
\end{lemma3}

Set
$$\mathbb F\triangleq \{X=(\zeta,U,\Gamma,R,H)\in\mathbb G \big| \zeta+Id+H\in G\}$$
and
$$\mathbb F_\alpha\triangleq \{X=(\zeta,U,\Gamma,R,H)\in\mathbb G \big| \zeta+Id+H\in G_\alpha\}.$$
Then we have the following lemma.

\begin{lemma3}
Let the map $\Phi: G\times\mathbb F\rightarrow\mathbb F,
(f,(\zeta,U,\Gamma,R,H))\mapsto (\bar{\zeta},\bar{U},\bar{\Gamma},\bar{R},\bar{H})$
be defined by
\begin{equation*}
\left\{\begin{array}{ll}
\bar{\zeta}=\zeta\circ f+f-Id,\\
\bar{U}=U\circ f,\\
\bar{\Gamma}=\Gamma\circ f,\\
\bar{R}=R\circ f,\\
\bar{H}=H\circ f.
\end{array}\right.
\end{equation*}
Then the map $\Phi$ defines a group action of G on $\mathbb F$.
\end{lemma3}

\begin{proof}
For any $f\in G$ and $X\in\mathbb F$, from the facts $(\zeta,U,\Gamma,H)\in (W^{1,\infty})^4$
and $f-Id\in W^{1,\infty}$, we see
$(\bar{\zeta},\bar{U},\bar{\Gamma},\bar{H})\in (W^{1,\infty})^4.$
According to the rule of chain,  one gets for almost everywhere $\xi\in \mathbb R$,
\begin{align} \label{3-lemma2-identities}
\bar{y}_\xi=y_\xi\circ ff_\xi,\ \ \bar{U}_\xi=U_\xi\circ ff_\xi,\ \ \bar{\Gamma}_\xi=\Gamma_\xi\circ ff_\xi,\ \ \bar{H}_\xi=H_\xi\circ ff_\xi.
\end{align}
Note that for any $f\in G$, there exists some large enough $\alpha>0$, such that $f\in G_\alpha$.
By Lemma 3.1, there exists some $c>0$ such that $\frac{1}{c}\leq f_\xi\leq c.$
Denote
\begin{align} \label{3-lemma2-denote}
\bar{X}\triangleq (\bar{\zeta},\bar{U},\bar{\Gamma},\bar{R},\bar{H})=(\zeta,U,\Gamma,R,H)\bullet f \triangleq X \bullet f.
\end{align}
By applying (\ref{3-lemma2-identities}), one can easily check that $\bar{X}$ satisfies (\ref{2-condition-2})-(\ref{2-condition-4}).

Now, we claim that $\bar{X}\in E.$ Indeed,
by recalling  $R\in L^\infty\cap L^2$ and $\frac{1}{c}\leq f_\xi\leq c,$
one obtains $\bar{R}\in L^\infty.$
Then by the change of variables, we have
\begin{eqnarray*}
\|\bar{R}\|^2_{L^2}
=\int_\mathbb R(R\circ f)^2(\xi)d\xi
=\int_\mathbb RR^2(x)f_\xi^{-1}(x)dx
\leq c\|R\|^2_{L^2}.
\end{eqnarray*}
So, $\bar{R}\in L^\infty\cap L^2.$
Similarly, by $U, \zeta_\xi\in L^2$, we get $\bar{U}, \bar{\zeta}_\xi\circ f \in L^2$.
Moreover,
$$\|\bar{\zeta}_\xi\|_{L^2}=\|\zeta_\xi\circ ff_\xi+f_\xi-1\|_{L^2}\leq c\|\bar{\zeta}_\xi\circ f\|_{L^2}+\|f_\xi-1\|_{L^2}<+\infty,$$
which along with the fact $\zeta\in W^{1,\infty}$ ensures $\zeta\in V.$
Likewise, $H\in V.$
Hence, $\bar{X}\in E.$
And thus, $\bar{X}\in \mathbb G.$

Next, notice that both $y+H$ and $f$ belong to the group $G$.
Then, $\bar{\zeta}+Id+\bar{H}=(y+H)\circ f\in G$.
So, $\bar{X}\in\mathbb F.$

On the other hand, in view of the definition of $\Phi$ and (\ref{3-lemma2-denote}), we infer that
\begin{eqnarray*}
&&X\bullet f_1\bullet f_2\\
&=&(\zeta\circ f_1+f_1-Id,U\circ f_1,\Gamma\circ f_1,R\circ f_1,H\circ f_1)\bullet f_2\\
&=&((\zeta\circ f_1+f_1-Id)\circ f_2+f_2-Id,U\circ f_1\circ f_2,\Gamma\circ f_1\circ f_2,
R\circ f_1\circ f_2,H\circ f_1\circ f_2)\\
&=&(\zeta\circ(f_1\circ f_2)+f_1\circ f_2-Id,U\circ(f_1\circ f_2),\Gamma\circ(f_1\circ f_2),
R\circ(f_1\circ f_2),H\circ(f_1\circ f_2))\\
&=&X\bullet(f_1\circ f_2),
\end{eqnarray*}
for any $X\in\mathbb F$ and $f_1,f_2\in G$.
Thus, the map $\Phi$ defines a group action of G on $\mathbb F$.
Therefore, we have proven the lemma.
\end{proof}

\begin{remark3}
By Lemma 3.2, we can consider the quotient space $\mathbb F/G$ of $\mathbb F$ with respect to the group action,
whose elements consist of $[X]$ defined by
$$[X]=\{X'\in\mathbb F \big| X' \sim X\},$$
where $X' \sim X$ means that there exists $f\in G$ such that $X'=X\bullet f.$
Moreover, for any $X\in\mathbb F,$ if we set
$\Pi(X)\triangleq X\bullet(\zeta+Id+H)^{-1}=X\bullet(y+H)^{-1}$,
then the mapping $\Pi$ is a projection, i.e. $\Pi\circ\Pi=\Pi$.
Hence, for any $X\in\mathbb F$ and $f\in G$, we have $T(X\bullet f)=T(X)$.
It follows that the mapping $[X]\rightarrow\Pi(X)$ is a bijection
from the quotient space $\mathbb F/G$ to $\mathbb F_0.$
\end{remark3}

Next, we turn to a property of the mapping $S_t$ in Theorem 2.1.
\begin{lemma3}
The mapping  $S_t$ is equi-variant. That is, for any $X\in\mathbb F$ and $f\in G,$ we have
$$S_t(X\bullet f)=S_t(X)\bullet f,$$
where $\bullet$ is defined in (\ref{3-lemma2-denote}).
\end{lemma3}

\begin{proof}
For any $X_0=(\zeta_0,U_0,\Gamma_0,R_0,H_0)\in\mathbb F$ and $f\in G,$
we denote
$\bar{X}_0=(\bar{\zeta}_0,\bar{U}_0,\bar{\Gamma}_0,\bar{R}_0,\bar{H}_0)=X_0\bullet f$, $X(t)=S_t(X_0)$
and $\bar{X}(t)=S_t(\bar{X}_0)$, respectively.

Now, we prove that $X(t)\bullet f$ satisfies (2.11) with the initial data $\bar{X}_0$.
%Indeed, by Theorem 2.1, the solution of (2.11) is unique, if $X(t)\bullet f$ and $\bar{X}(t)$ satisfy the same system of differential equation with the same initial data,  then they are equal.
For this, denote
$$\hat{X}(t)=(\hat{\zeta}(t),\hat{U}(t),\hat{\Gamma}(t),\hat{R}(t),\hat{H}(t))=X(t)\bullet f.$$
Since $X(t)$ is the solution to system (\ref{2-new-system}), it follows that
\begin{align*}
\hat{\zeta}(t)_t=(\zeta(t)\circ f+f-Id)_t=\zeta_t\circ f=U(t)\circ f=\hat{U}(t),
\end{align*}
and
\begin{eqnarray*}
\hat{U}(t)_t
&=&(U(t)\circ f)_t\\
&=& U_t(t)\circ f\\
&=& -G_1(X)\circ f\\
&=&\frac{1}{4}\int_\mathbb R sgn(f(\xi)-\eta)\exp\left(-|y(f(\xi))-y(\eta)|\right)[U^2y_\xi-2R^2y_\xi+H_\xi](\eta)d\eta,
\end{eqnarray*}
%\begin{align*}
%\hat{U}(t)_t=&(U(t)\circ f)_t=U_t(t)\circ f=-G_1(X)\circ f\\
%=&\frac{1}{4}\int_\mathbb Rsgn(f(\xi)-\eta)exp\left(-|y(f(\xi))-y(\eta)|\right)[U^2y_\xi-2R^2y_\xi+H_\xi](\eta)d\eta,
%\end{align*}
where $G_1$ is defined in Lemma 2.1.
Thanks to $\hat{y}(t)=y(t)\circ f$ and $\hat{H}(t)=H(t)\circ f,$ together with the chain rule,
we get
$$\hat{y}_\xi(t,\xi)=y_\xi(f(\xi))f_\xi(\xi)\quad \textrm{and}\quad \hat{H}_\xi(t,\xi)=H_\xi(f(\xi))f_\xi(\xi).$$
Applying the change of variables, and noting that $f$ is increasing, one infers
\begin{align*}
&\hat{U}_t(t)\\
=&\frac{1}{4}\int_\mathbb R sgn(f(\xi)-f(\eta))
\exp\left(-|y(f(\xi))-y(f(\eta))|\right)[U^2y_\xi-2R^2y_\xi+H_\xi](f(\eta))f_\xi(\eta)d\eta\nonumber\\
=&\frac{1}{4}\int_\mathbb R sgn(\xi-\eta)
\exp\left(-|\hat{y}(\xi)-\hat{y}(\eta)|\right)[\hat{U}^2\hat{y}_\xi-2\hat{R}^2\hat{y}_\xi+\hat{H}_\xi](\eta)d\eta\nonumber\\
=&-G_1(\hat X).
\end{align*}
Since the other terms in system (\ref{2-new-system}) can be treated in a similar way,
it follows that $\hat{X}(t)$ is a solution to system (\ref{2-new-system}) with the initial $\hat{X}(0)=\bar{X}_0$.
Then Theorem 2.1 ensures
$$ S_t(X_0)\bullet f=X(t)\bullet f=\hat{X}(t)=\bar{X}(t)=S_t(\bar{X}_0)=S_t(X_0\bullet f),$$
which completes the proof of the lemma.
\end{proof}

\begin{remark3}
By Lemma 3.3 and Remark 3.1, we have
\begin{align*}
\Pi\circ S_t\circ\Pi=\Pi\circ S_t.
\end{align*}
If we define the semigroup $\bar{S}_t$ on $\mathbb F_0$ as
\begin{align} \label{3-remark-2}
\bar{S}_t=\Pi\circ S_t,
\end{align}
Then by Theorem 2.1, $S_t$ is a continuous semigroup. So is $\bar{S}_t$.
\end{remark3}

To obtain the Lipschitz continuity of the solutions to system (\ref{2-new-system}),
we need to introduce the distances on $\mathbb F$ and $\mathbb F_0,$ respectively.
Precisely, for any $X_\alpha, X_\beta\in\mathbb F$,
we define
\begin{equation}\label{3-distance-J}
J(X_\alpha,X_\beta)\triangleq \inf_{f,g\in G}\|X_\alpha\bullet f-X_\beta\|_E+\|X_\alpha-X_\beta\bullet g\|_E.
\end{equation}
While for any $X_\alpha,X_\beta\in\mathbb F_0$, we set
\begin{equation}\label{3-distance-d}
d(X_\alpha,X_\beta)\triangleq \inf\sum^{n=N}_{n=1}J(X_{n-1},X_n),
\end{equation}
where the infimum is taken over all finite sequences $\{X_n\}_{n=0}^N\subset\mathbb F_0$
satisfying $X_0=X_\alpha$ and $X_N=X_\beta.$

Similar to the arguments in \cite{G-H-Rdcds}, we can readily get the following result.
\begin{lemma3}
(i) For $X_\alpha,X_\beta\in\mathbb F$ and $f\in G_k,$
$$\|X_\alpha\bullet f-X_\beta\bullet f\|_E\leq C\|X_\alpha-X_\beta\|_E
\ \ \textrm{and} \ \
J(X_\alpha\bullet f,X_\beta)\leq CJ(X_\alpha,X_\beta),$$
where $C=C(k)>0$.

(ii) For any $X_\alpha,X_\beta\in\mathbb F_0$, we have
\begin{align*}
\frac{1}{2}\|X_\alpha-X_\beta\|_{L^\infty}\leq d(X_\alpha,X_\beta)\leq2\|X_\alpha-X_\beta\|_E.
\end{align*}
\end{lemma3}

For any $M>0,$ set
$$\mathbb F^M \triangleq \{X=(\zeta,U,\Gamma,R,H)\in\mathbb F \big| \|H\|_{L^\infty}\leq M\}$$
and
$$\mathbb F_0^M \triangleq \mathbb F_0\cap\mathbb F^M.$$
\begin{remark3}
$\mathbb F_0^M$ and  $\mathbb F_0\cap B_M$ are equivalent in the sense that
$$\mathbb F_0\cap B_M\subset\mathbb F_0^M\subset\mathbb F_0\cap B_{M'},$$
where $B_M \triangleq \{X\in E|\|X\|_E\leq M\}$, and $M'$ depends only on $M$.
\end{remark3}

Let us define a distance on $\mathbb F_0^M$ as follows:
\begin{equation} \label{3-distance-dM}
d^M(X_\alpha,X_\beta)\triangleq \inf\sum^{n=N}_{n=1}J(X_{n-1},X_n),
\end{equation}
where the infimum is taken over all finite sequences $\{X_n\}_{n=0}^N\subset\mathbb F_0^M$
satisfying $X_0=X_\alpha$ and $X_N=X_\beta.$

Now, we are in a position to state our main theorem in this section.

\begin{theorem3}
Given $T>0$ and $M>0,$ there exists a $C=C(T,M)>0$ such that
for any $X_\alpha,X_\beta\in\mathbb F_0^M$ and $t\in[0,T]$,
\begin{align*}
d^M(\bar{S}_t(X_\alpha),\bar{S}_t(X_\beta))\leq C d^M(X_\alpha,X_\beta),
\end{align*}
where $\bar{S}_t$ is defined by (\ref{3-remark-2}).
\end{theorem3}

\begin{proof}
By the definitions of $d^M$ and $J$, for arbitrary $\epsilon\in(0,1)$,
there exist sequences $\{X_n\}_{n=0}^N\subset\mathbb F_0, \{\tilde{f}_n\}_{n=1}^{N},\{f_n\}_{n=0}^{N-1}\subset G$
with $X_0=X_\alpha, X_N=X_\beta$ such that
\begin{align} \label{3-thm-1}
\sum_{n=1}^N(\|X_{n-1}\bullet f_{n-1}-X_n\|_E+\|X_{n-1}-X_n\bullet\tilde{f}_n\|_E)\leq d^M(X_\alpha,X_\beta)+\epsilon.
\end{align}
Denote
$$X_n^t=S_t(X_n),\,\ g_n^t=\zeta_n^t+Id+H_n^t, \,\ \bar{X}_n^t=\bar{S}_tX_n=\Pi(X_n^t)=X_n^t\bullet(g_n^t)^{-1}.$$
According to Lemma 3.3 and the similar arguments of Lemma 2.5 in \cite{G-H-Rdcds},
for $k,T>0$ and $X\in\mathbb F_k$, there exists a $k'$ depending only on $k,T$ and $\|X\|_E$,
such that $S_t(X)\in\mathbb F_{k'}$. Then $g_n^t\in G_k$ for some $k=k(M,T).$
Thanks to $X_n\in\mathbb F^M$, we get $\bar{X}_n^t\in\mathbb F_0^M.$ \\
For $f_{n-1}^t=g_{n-1}^t\bullet f_{n-1}\bullet(g_n^t)^{-1},$ by using Lemma 3.4, one deduces that
\begin{align*}
\|\bar{X}_{n-1}^t\bullet f_{n-1}^t-\bar{X}_n^t\|_E&=\|X_{n-1}^t\bullet(g_{n-1}^t)^{-1}\bullet f_{n-1}^t-X_n^t\bullet(g_n^t)^{-1}\|_E\\
&\leq C \|X_{n-1}^t\bullet(g_{n-1}^t)^{-1}\bullet f_{n-1}^t\bullet(g_n^t)-X_n^t\|_E\\
&=C\|X_{n-1}^t\bullet f_{n-1}-X_n^t\|_E\\
&=C\|S_t(X_{n-1}\bullet f_{n-1})-S_t(X_n)\|_E.
\end{align*}

On the other hand, by Remark 3.3, there is some $M'>0$ such that
$\|X\|_E\leq M'$ for any $X\in\mathbb F_0^M.$
Then (\ref{3-thm-1}) and the triangle inequality imply that
\begin{align*}
\|X_{n-1}\bullet f_{n-1}\|_E&\leq\|X_{n-1}\bullet f_{n-1}-X_n\|_E+\|X_n\|_E\\
&\leq d^M(X_\alpha,X_\beta)+\epsilon+M'\\
&\leq2\|X_\alpha-X_\beta\|_E+\epsilon+M'\\
&\leq5M'+1,
\end{align*}
which along with Theorem 2.1 yields
$$\|\bar{X}_{n-1}^t\bullet f_{n-1}^t-\bar{X}_n^t\|_E\leq C\|X_{n-1}t\bullet f_{n-1}-X_n\|_E,$$
where $C=C(T,M)>0$.

Likewise, for $\tilde{f}_n^t=g_n^t\bullet\tilde{f}_n\bullet(g_{n-1}^t)^{-1}$, there holds
$$\|\bar{X}_{n-1}^t-\bar{X}_n^t\bullet\tilde{f}_n^t\|_E\leq C_M\|X_{n-1}t-X_n\bullet \tilde{f}_{n}\|_E.$$
Thus, by using (\ref{3-distance-dM}) and (\ref{3-thm-1}), one has
\begin{align*}
&d^M(\bar{S}_t(X_\alpha),\bar{S}_t(X_\beta))\\
\leq&\sum_{n=1}^N(\|\bar{X}_{n-1}^t\bullet f_{n-1}^t-\bar{X}_n^t\|_E+\|\bar{X}_{n-1}^t-\bar{X}_n^t\bullet\tilde{f}_n^t\|_E)\\
\leq& C\sum_{n=1}^N(\|X_{n-1}t\bullet f_{n-1}-X_n\|_E+\|X_{n-1}t-X_n\bullet \tilde{f}_{n}\|_E)\\
\leq& C(d^M(X_\alpha,X_\beta)+\epsilon),
\end{align*}
which leads to the desired result.
\end{proof}

\section{Global weak solutions to the system (\ref{1-M2CH})}
\newtheorem {remark4}{Remark}[section]
\newtheorem{theorem4}{Theorem}[section]
\newtheorem{lemma4}{Lemma}[section]
\newtheorem{definition4}{Definition}[section]

In this section, we shall show that the original system (\ref{1-M2CH}) has global weak solutions as
the initial data $(u_0,\gamma_0)\in H^1\times (H^1\cap W^{1,\infty})$.
To achieve this, we first introduce some definitions as follows.

\begin{definition4}
Given initial data $z_{0}=(u_{0}, \gamma_{0})\in H^1\times(H^{1}\cap W^{1,\infty}).$
If $z(t,x)=(u(t,x),\gamma(t,x))\in L^{\infty}((0,\infty); H^{1}\times H^{1})$
satisfies system (\ref{1-M2CH}) and
$z(t,\cdot)\rightarrow z_0$ as $t\rightarrow 0^+$ in the sense of distribution,
then $z=(u,\gamma)$ is called a global weak solution to system (\ref{1-M2CH}).
Moreover, if
$$h(t)\triangleq \int_\mathbb R (u^2+u_x^2+\gamma^2+\gamma_x^2)(t,x)dx=h(0)$$
holds for almost all $t>0,$ then $z$ is called a global conservative weak solution.
\end{definition4}

%We give a set in Eulerian coordinates.

\begin{definition4}
The set $\mathbb D$ consists of all  $(u,\gamma,\mu)$ such that \\
(i) $(u,\gamma)\in H^1\times(H^1\cap W^{1,\infty}),$ \\
(ii) $\mu$ is a positive Radon measure whose absolutely continuous part $\mu_{ac}$ satisfies
\begin{align*} \label{4-def-2}
\mu_{ac}=u^2+u_x^2+\gamma^2+\gamma_x^2.
\end{align*}
\end{definition4}

\begin{definition4} \cite{H-R1}
For any $(u,\gamma,\mu)\in\mathbb D,$ define the map
$L: \mathbb D\rightarrow \mathbb F_0, (u,\gamma,\mu)\mapsto X\triangleq(\zeta,U,\Gamma,R,H)$ by
\begin{equation*} \label{4-def-3}
\left\{\begin{array}{ll}
\zeta(\xi)=y(\xi)-\xi=\sup\{y \big| \mu(-\infty,y)+y<\xi\}-\xi,\\
U(\xi)=u\circ y(\xi),\\
\Gamma(\xi)=\gamma\circ y(\xi),\\
R(\xi) =\gamma_x\circ y(\xi),\\
H(\xi)=-\zeta(\xi).
\end{array}\right.
\end{equation*}
\end{definition4}

As a result, for any initial data $(u_0,\gamma_0,\mu_0)\in\mathbb D,$
we can construct a solution to system (\ref{2-new-system}) in $\mathbb F$
with initial data $X_0=L(u_0,\gamma_0,\mu_0)\in\mathbb F_0.$
Thus, the global weak solutions to system (\ref{2-new-system})
yields a global conservative weak solution to system (\ref{1-M2CH}) in the original variables,
which is the goal of this section.

\begin{definition4} \cite{H-R1}
The mapping $W: \mathbb F_0 \rightarrow \mathbb D, X\triangleq(\zeta,U,\Gamma,R,H)\mapsto (u,\gamma,\mu)$
is defined by
$$u(x)=U(\xi),\ \gamma(x)=\Gamma(\xi)$$
for any $\xi$, such that $x=y(\xi)$ and
$$\mu(B)=\int_{\{x\in y^{-1}(B)\}} H_\xi(x)dx$$ for any Borel set $B$.
\end{definition4}

\begin{remark4}
(i) \begin{align} \label{4-Id}
L\circ W=Id_{\mathbb F_0}\ \ \ \ \textrm{and}\ \ \ W\circ L=Id_{\mathbb D}.
\end{align}
(ii) The distance $J$ in (\ref{3-distance-J}) can be viewed as a map from $\mathbb F/G$ to $\mathbb D$, i.e.
\begin{align} \label{4-W-identity}
W(X_1)=W(X),\  \forall\, X_1\in[X].
\end{align}
\end{remark4}

Set
\begin{align} \label{4-Tt}
T_t\triangleq W\circ \bar{S}_t \circ L: \mathbb D \rightarrow \mathbb D,
\end{align}
where $\bar{S}_t$ is defined in (\ref{3-lemma2-denote}).
Then we can prove our main theorem in this section.

\begin{theorem4}
Given any initial data $(u_0,\gamma_0,\mu_0)\in\mathbb D$,
denote $(u,\gamma,\mu)\triangleq T_t(u_0,\gamma_0,\mu_0)$.
Then $z=(u,\gamma)$ is a global conservative weak solution to system (\ref{1-M2CH}).
\end{theorem4}

\begin{proof}
Since Lemma 3.3, it follows from (\ref{4-W-identity}) that $T_t=W\circ S_t\circ L$.
Let
$$X(t)=(\zeta(t),U(t),\Gamma(t),R(t),H(t))=S_t(L(X_0)).$$
Then $X(t)$ is the solution to system (\ref{2-new-system}) with initial data $L(X_0).$
So, for any smooth function $\phi$ with the compact support in $\mathbb R_+\times\mathbb R,$  we have
\begin{align*}
 &-\int_{\mathbb R_+\times\mathbb R}u(t,x)\phi_t(t,x)dxdt\\
 =&-\int_{\mathbb R_+\times\mathbb R}u(t,y(t,\xi))\phi_t(t,y(t,\xi))y_\xi(t,\xi)d\xi dt\\
 =&-\int_{\mathbb R_+\times\mathbb R}U(t,\xi)[(\phi(t,y(t,\xi)))_t-\phi_x(t,y(t,\xi))y_t(t,\xi)]y_\xi(t,\xi)d\xi dt\\
 =&-\int_{\mathbb R_+\times\mathbb R}(U(t,\xi)(\phi(t,y(t,\xi)))_ty_\xi(t,\xi)-U^2(t,\xi)\phi_x(t,y(t,\xi))y_\xi(t,\xi))d\xi dt\\
 =&\int_{\mathbb R}U(0,\xi)\phi(0,y(0,\xi))y_\xi(0,\xi)d\xi+\int_{\mathbb R_+\times\mathbb R}(U(t,\xi)y_\xi(t,\xi))_t\phi(t,y(t,\xi))d\xi dt\\
 &+\int_{\mathbb R_+\times\mathbb R}U^2(t,\xi)(\phi(t,y(t,\xi))_\xi d\xi dt\\
 =&\int_{\mathbb R}U(0,\xi)\phi(0,y(0,\xi))y_\xi(0,\xi)d\xi\\&+\int_{\mathbb R_+\times\mathbb R}(U_t(t,\xi)y_\xi(t,\xi)+U(t,\xi)y_{t\xi}(t,\xi))\phi(t,y(t,\xi))d\xi dt\\&-2\int_{\mathbb R_+\times\mathbb R}U(t,\xi)U_\xi(t,\xi)\phi(t,y(t,\xi)d\xi dt\\
 =&\int_{\mathbb R}U(0,\xi)\phi(0,y(0,\xi))y_\xi(0,\xi)d\xi\\&+\int_{\mathbb R_+\times\mathbb R}(U_t(t,\xi)y_\xi(t,\xi)+U(t,\xi)U_\xi(t,\xi))\phi(t,y(t,\xi))d\xi dt\\&-2\int_{\mathbb R_+\times\mathbb R}U(t,\xi)U_\xi(t,\xi)\phi(t,y(t,\xi)d\xi dt\\
 =&\int_{\mathbb R}U(0,\xi)\phi(0,y(0,\xi))y_\xi(0,\xi)d\xi\\
 &-\int_{\mathbb R_+\times\mathbb R}(Q_1(t,\xi)y_\xi(t,\xi)+U(t,\xi)U_\xi(t,\xi))\phi(t,y(t,\xi))d\xi dt,
 \end{align*}
and
\begin{align*}
 &-\int_{\mathbb R_+\times\mathbb R}(\frac{1}{2}u^2(t,x)+P_1(t,x))\phi_x(t,x)dxdt\\
 =&\int_{\mathbb R_+\times\mathbb R}(u(t,y(t,\xi))u_x(y,y(t,\xi))+P_{1,x}(t,y(t,\xi))\phi(t,y(t,\xi))y_\xi(t,\xi)d\xi dt\\
 =&\int_{\mathbb R_+\times\mathbb R}(U(t,\xi)U_\xi(t,\xi)+P_{1,\xi}(t,\xi))\phi(t,y(t,\xi))d\xi dt\\
 =&\int_{\mathbb R_+\times\mathbb R}(U(t,\xi)U_\xi(t,\xi)+Q_1(t,\xi)y_\xi(t,\xi))\phi(t,y(t,\xi))d\xi dt.
\end{align*}
Hence, we get
\begin{align} \label{4-thm-1}
 &\int_{\mathbb R_+\times\mathbb R}(-u(t,x)\phi_t(t,x)-(\frac{1}{2}u^2(t,x)+P_1(t,x))\phi_x(t,x))dxdt\\
 =&\int_{\mathbb R}U(0,\xi)\phi(0,y(0,\xi))y_\xi(0,\xi)d\xi\nonumber\\
 =&\int_\mathbb Ru_0(x)\phi(0,x)dx.\nonumber
\end{align}

Next, we show that
\begin{align} \label{4-thm-2}
 P_1(t,x)-P_{1,xx}(t,x)=u^2+\frac{1}{2}u_x^2+\frac{1}{2}\gamma^2-\frac{1}{2}\gamma_x^2
\end{align}
holds in the sense of distribution.
Indeed, by using (\ref{2-equality-two-systems}), (\ref{2-condition-3}) and Lemma 2.4,  one infers that
\begin{align*}
 &\int_{\mathbb R_+\times\mathbb R}P_{1,x}(t,x))\phi_x(t,x)dxdt\\
 \nonumber=&\int_{\mathbb R_+\times\mathbb R}P_{1,x}(t,y(t,\xi))\phi_x(t,y(t,\xi))y_\xi(t,\xi)d\xi dt\\
 \nonumber=&\int_{\mathbb R_+\times\mathbb R}P_{1,\xi}(t,\xi)\phi_x(t,y(t,\xi))d\xi dt\\
 \nonumber=&\int_{\mathbb R_+\times\mathbb R}Q_1(t,\xi)y_\xi(t,\xi)\phi_x(t,y(t,\xi))d\xi dt\\
 \nonumber =&\int_{\mathbb R_+\times\mathbb R}Q_1(t,\xi)\phi_\xi(t,y(t,\xi))d\xi dt\\
 \nonumber=&-\int_{\mathbb R_+\times\mathbb R}Q_{1,\xi}(t,\xi)\phi(t,y(t,\xi))d\xi dt\\
 \nonumber=&\int_{\mathbb R_+\times\mathbb R}[\frac{1}{2}U^2y_\xi-R^2y_\xi-P_1(t,\xi)y_\xi+\frac{1}{2}H_\xi]\phi(t,y(t,\xi))d\xi dt\\
 \nonumber=&\int_{\{(t,\xi)\in\mathbb R_+\times\mathbb R|y_\xi>0\}}\left((u^2+\frac{1}{2}u_x^2+\frac{1}{2}\gamma^2-\frac{1}{2}\gamma_x^2)\circ y(t,\xi)y_\xi(t,\xi)-P_1y_\xi\right)\phi(t,y(t,\xi))d\xi dt\\
 \nonumber=&\int_{\mathbb R_+\times\mathbb R}\left((u^2+\frac{1}{2}u_x^2+\frac{1}{2}\gamma^2-\frac{1}{2}\gamma_x^2)\circ y(t,\xi)-P_1(t,\xi)\right)\phi(t,y(t,\xi))y_\xi d\xi dt\\
 \nonumber=&\int_{\mathbb R_+\times\mathbb R}((u^2+\frac{1}{2}u_x^2+\frac{1}{2}\gamma^2-\frac{1}{2}\gamma_x^2)-P_1(t,x))\phi(t,x)dxdt,
\end{align*}
which implies that (\ref{4-thm-2}) holds in the sense of distributions.
From (\ref{4-thm-1}) and (\ref{4-thm-2}), one can easily check that
$(u,\gamma)$ satisfies system (\ref{1-M2CH}) in the sense of distribution on $\mathbb R_+\times\mathbb R.$

On the one hand, by Definition 4.4 and Theorem 2.1, we have
\begin{align} \label{4-thm-3}
 \mu(t)(\mathbb R)=\int_\mathbb RH_\xi(\xi)d\xi=H(t,+\infty)=H(0,+\infty).
\end{align}
Note that Lemma 2.4 implies that there exists $K\subset\mathbb R_+$ with $meas(K^c)=0,$
such that $y_\xi(t,\xi)>0$ for any $t\in K$ and almost every $\xi\in\mathbb R.$
Given $t\in K$ and any Borel set $B,$ thanks to (\ref{2-condition-3}) again, we obtain
\begin{align*}
 \mu(t)(B)=&\int_{y^{-1}(B)}H_\xi d\xi\\
 =&\int_{y^{-1}(B)}(U^2+\Gamma^2+R^2+\frac{U_\xi^2}{y_\xi})y_\xi d\xi\\
 =&\int_B(u^2+u_x^2+\gamma^2+\gamma_x^2)(t,x)dx,
\end{align*}
which along with (\ref{4-thm-3}) yields that
$$h(t)=\int_\mathbb R(u^2+u_x^2+\gamma^2+\gamma_x^2)(t,x)dx=h(0)$$
holds for almost all $t\in\mathbb R_+.$
Therefore, we have completed the proof of the theorem.
 \end{proof}

\section{Lipschitz continuity of the weak solution}
\newtheorem {remark5}{Remark}[section]
\newtheorem{theorem5}{Theorem}[section]
\newtheorem{lemma5}{Lemma}[section]

In this section, we consider the Lipschitz continuity and uniqueness of the weak solution.
To this end, let us first define the metric $d_\mathbb D:\mathbb D\times\mathbb D\rightarrow\mathbb R_+$ by
\begin{align*}
d_\mathbb D((u_1,\gamma_1,\mu_1),(u_2,\gamma_2,\mu_2))=d(L(u_1,\gamma_1,\mu_1),L(u_2,\gamma_2,\mu_2)),
\end{align*}
where $d$ and $L$ are defined in (\ref{3-distance-d}) and Definition 4.3, respectively.

Moreover, another metric $d_{\mathbb D^M}:\mathbb D^M\times\mathbb D^M\rightarrow\mathbb R_+$ is defined by
\begin{align}  \label{5-def-2}
d_{\mathbb D^M}((u_1,\gamma_1,\mu_1),(u_2,\gamma_2,\mu_2))=d_M(L(u_1,\gamma_1,\mu_1),L(u_2,\gamma_2,\mu_2)),
\end{align}
where $\mathbb D^{M}=\{(u,\gamma,\mu)\in\mathbb D\mid\mu(\mathbb R)\leq M\}$ for some $M>0,$
and the metric $d_M$ is defined in (\ref{3-distance-dM}).
Note that if $(u,\rho,\mu)\in\mathbb D^M$, then $L(u,\rho,\mu)\in\mathbb F_0^M.$
Thus the metric in (\ref{5-def-2}) is well defined.

The main result of this section is stated as follows:

\begin{theorem5}
The semigroup $(T_t,d_\mathbb D)$ is a continuous semigroup on $\mathbb D$
and satisfies the following Lipschitz continuity:
for any  $M,\, T>0$, there exists some $C=C(M,T)>0$ such that
for any $(u_1,\gamma_1,\mu_1),(u_2,\gamma_2,\mu_2)\in\mathbb D^M$ and $t\in[0,T]$,
$$d_{\mathbb D^M}(T_t(u_1,\gamma_1,\mu_1),T_t(u_2,\gamma_2,\mu_2))
\leq C d_{\mathbb D^M}((u_1,\gamma_1,\mu_1),(u_2,\gamma_2,\mu_2)).$$
\end{theorem5}

\begin{proof}
We first prove that $T_t$ is a semigroup.
Indeed, since $\bar{S}_t$ is a mapping from $\mathbb F_0$ to $\mathbb F_0$,
it follows from (\ref{4-Id}) that
$$T_t\circ T_{t'}=W\circ\bar{S}_tL\circ W\circ\bar{S}_{t'}\circ L=W\circ\bar{S}_t\circ\bar{S}_{t'}\circ L=W\circ\bar{S}_{t+t'}\circ L=T_{t+t'}.$$
On the other hand, thanks to Theorem 3.1 and (\ref{4-Id}), we have
\begin{align*}
&d_{\mathbb D^M}(T_t(u_1,\gamma_1,\mu_1),T_t(u_2,\gamma_2,\mu_2))\\
=&d^M(L(T_t(u_1,\gamma_1,\mu_1)),L(T_t(u_2,\gamma_2,\mu_2)))\\
=&d^M(\bar{S}_t\circ L(u_1,\gamma_1,\mu_1),\bar{S}_t\circ L(u_2,\gamma_2,\mu_2))\\
\leq& C d^M((L(u_1,\gamma_1,\mu_1)),(L(u_2,\gamma_2,\mu_2)))\\
=&C d_{\mathbb D^M}((u_1,\gamma_1,\mu_1),(u_2,\gamma_2,\mu_2)),
\end{align*}
which completes the proof of the theorem.
\end{proof}

Similar to the arguments in \cite{H-R1}, we can readily prove the following lemma.

\begin{lemma5}
Let $(u,\gamma)$ be the weak solution to system (\ref{1-M2CH}).
Then the mapping
$$(u,\gamma)\mapsto (u,\gamma,(u^2+u_x^2+\gamma^2+\gamma^2_x)dx)$$
is continuous from $H^1\times H^1$ into $\mathbb D.$
\end{lemma5}

By applying Lemma 5.1 and Theorem 5.1, one can immediately get the following uniqueness result.

\begin{theorem5}
For any initial data $z_0=(u_0,\gamma_0)\in H^1\times (H^1\cap W^{1,\infty}),$
the system (\ref{1-M2CH}) has a unique global conservative weak solution
$z=(u,\gamma)\in H^1\times H^1$ with the following form
$$(u,\gamma,(u^2+u_x^2+\rho^2)dx)=T_t(z_0),$$
where the mapping $T_t$ is defined by (\ref{4-Tt}).
\end{theorem5}

\bigskip

\textbf{Acknowledgments}
The authors thank the referees for their valuable comments and suggestions.
Guan was partially supported by NSFY (No.11201494).
Yan was partially supported by NNSFC (No.11501226).
Wei was partially supported by CNNSF (No.11101095) and the High-Level Talents Project of Guangdong Province (No.2014011).

\end{document}